\documentclass[11pt,reqno]{amsart}
\usepackage{amsmath}
\usepackage{amsfonts}
\usepackage{amssymb}

\newtheorem{theorem}{Theorem}
\theoremstyle{plain}
\newtheorem{lemma}{Lemma}
\newtheorem{corollary}{Corollary}

\theoremstyle{definition}

\newtheorem*{acknowledgements}{Acknowledgements}
\numberwithin{equation}{section}

\newcommand{\thmref}[1]{Theorem~\ref{#1}}

\newcommand{\lemref}[1]{Lemma~\ref{#1}}

\newcommand{\by}[1]{\mathrm{by\ (\ref{#1})}}

\newcommand{\Mlt}{\mathrm{Mlt}}
\newcommand{\Nuc}{\mathrm{Nuc}}

\begin{document}
\title{Every diassociative A-loop is Moufang}
\author{Michael K. Kinyon}
\address{Department of Mathematics \& Computer Science\\
Indiana University\\
South Bend, IN\ 46634 USA}
\email{mkinyon@iusb.edu}
\urladdr{http://www.iusb.edu/\symbol{126}mkinyon}
\author{Kenneth Kunen}
\thanks{The second author's work was
partly supported by NSF Grant DMS-9704520.} 
\address{Department of Mathematics\\
University of Wisconsin\\
Madison, WI 57306 USA}
\email{kunen@math.wisc.edu}
\urladdr{http://www.math.wisc.edu/\symbol{126}kunen}
\author{J.D. Phillips}
\address{Department of Mathematics\\
Saint Mary's College of CA\\
Moraga, CA 94575\ USA}
\email{phillips@stmarys-ca.edu}
\date{\today}

\begin{abstract}
An A-loop is a loop in which every inner mapping is an automorphism.
We settle a problem which had been open since 1956 by showing that
every diassociative A-loop is Moufang.
\end{abstract}\maketitle

\section{Introduction}
\label{sec:intro}

A \emph{loop} $(L,\cdot)$ consists of a nonempty set $L$
with a binary operation $\cdot$ on $L$ such that (i) given
$a,b\in L$, the equations $ax=b$ and $ya=b$ each have unique
solutions $x,y\in L$, and (ii) there exists an identity element
$1\in L$ satisfying $1x=x1=x$ for all $x\in L$. As usual, we
abbreviate the binary operation by juxtaposition. Two varieties
of loops which have been widely discussed in the literature
are the \textit{Moufang loops} and the \textit{A-loops}.

A \emph{Moufang loop} is a loop satisfying the identity
\begin{equation}\label{eq:moufang}
x(y\cdot xz)=(xy\cdot x)z.
\end{equation}
These were introduced by R. Moufang in 1934 \cite{moufang},
and are discussed in detail in the texts by Bruck \cite{bruck}
and Pflugfelder \cite{pflugfelder}. By Moufang's Theorem
(\cite{bruck}, VII.4; \cite{pflugfelder}, IV.2.9), every Moufang
loop is diassociative; that is, the subloop $\langle x,y \rangle$
generated by any pair of elements is a group.

For $x\in L$, the left and right translations by $x$ are defined
by $yL(x)=xy$ and $yR(x)=yx$, respectively. The
\emph{multiplication group} of $L$ is the permutation group
$\Mlt(L) = \langle R(x),L(x):x\in L\rangle$ generated by all left and
right translations. The \emph{inner mapping group} is the subgroup
$\Mlt_1(L)$ fixing $1$. If $L$ is a group, then $\Mlt_1(L)$ is the
group of inner automorphisms of $L$.

In 1956, R.H. Bruck and L.J. Paige \cite{bp} defined an
\emph{A-loop} to be a loop in which every inner mapping is
an automorphism. Many of the basic theorems about A-loops
are contained in \cite{bp}; for example, A-loops are always
power associative (every $\langle x \rangle$ is a group), but
not necessarily diassociative. In the same paper, Bruck and
Paige included a detailed study of the diassociative A-loops,
pointing out that these satisfy ``many of the properties of
Moufang loops''. In hindsight, this is not surprising, since,
as we will show:

\begin{theorem}\label{thm:main}
Every diassociative A-loop is a Moufang loop.
\end{theorem}

For commutative loops, this was proved in 1958 by J.M. Osborn
\cite{osborn}. Conversely, every commutative Moufang loop is an
A-loop (see Bruck \cite{bruck}, Lemma VII.3.3). However, not
all Moufang loops are A-loops; \cite{bp, phillips}, together with
the results of the present paper, provide a simple description
of the diassociative A-loops as a sub-variety of the Moufang loops
(see Corollary \ref{cor:moufA}). Further work
on A-loops is contained in \cite{drapal, fps}.

By our \thmref{thm:main}, we have:

\begin{corollary}
\label{cor:equiv}
For an A-loop, the following are equivalent:
\begin{enumerate}
\item $L$ \emph{has the} inverse property, \emph{i.e.,} $x^{-1}(xy)=y$
\emph{and} $(xy)y^{-1}=x$ \emph{for all} $x,y\in L$;

\item $L$ \emph{has the} alternative property, \emph{i.e.,} $x(xy)=x^{2}y$
\emph{and} $(xy)y=xy^{2}$ \emph{for all} $x,y\in L$;

\item $L$ \emph{is diassociative.}

\item $L$ \emph{is a Moufang loop}.
\end{enumerate}
\end{corollary}

The equivalence of the first three items is from
Bruck and Paige \cite{bp}, Theorem~3.1. (One may begin
with even weaker hypotheses, but we will not pursue
this here.) 

In any loop, the inner mapping group $\Mlt_1(L)$ is generated by the
left, right, and middle inner mappings defined, respectively, by:
\begin{align*}
L(x,y) &= L(x)L(y)L(yx)^{-1}\\
R(x,y) &= R(x)R(y)R(xy)^{-1}\\
T(x) &= R(x)L(x)^{-1}
\end{align*}
(\cite{bruck}, IV.1, \cite{pflugfelder}, I.5.2). Bruck and Paige
(\cite{bp}, (3.42)) showed that diassociative A-loops satisfy:
\begin{equation}\label{eq:bp(3.42)}
\left[ R(x,y)R(y,x)\right]^{-1}T(x)T(y)=T(xy)
\end{equation}
Furthermore, they showed (see Corollary on p.~315) that for
Moufang A-loops, the map $T : L \to \Mlt_1(L)$ (where $x\mapsto T(x)$)
is a homomorphism (i.e., $T(x)T(y)=T(xy)$, so $R(x,y) = R(y,x)^{-1}$).
Not surprisingly, one of our key lemmas will be:

\begin{lemma}\label{lem:T-homom}
If $L$ is a diassociative A-loop, then $T : L \to \Mlt_1(L)$ is a
\mbox{homomorphism.}
\end{lemma}

The nucleus, $\Nuc(L)$, of an inverse property loop
$L$ is the normal subloop of all elements that associate with
all pairs of elements from $L$, i.e.,
$\Nuc (L)= \{x\in L:(xy)z=x(yz)$ for all $y,z\in L\}$.
By results already in the literature, we have the following
corollary to Theorem \ref{thm:main}:

\begin{corollary} \label{cor:moufA}
$L$ is a diassociative A-loop if and only if $L$ is Moufang and 
$L/\Nuc(L)$ is a commutative loop of exponent three.
\end{corollary}

\begin{proof}
In any Moufang loop, each $T(x)$ is a pseudo-automorphism with
companion $x^{-3}$, and each $R(x,y) = L(x^{-1},y^{-1})$ is
a pseudo-automorphism with companion the commutator $(x,y)$
(\cite{bruck}, Lemma VII.2.2). In general, if $c$ is a companion
of the pseudo-automorphism $\varphi$, then $c$ is in the nucleus
iff $\varphi$ is an automorphism. Thus all cubes and commutators
are in the nucleus iff all inner mappings are automorphisms.
\end{proof}

Every Moufang A-loop is an $M_4$ loop in the terminology of
Pflugfelder \cite{hala, pflugfelder}; that is, it satisfies
the identity $(xy)(zx^4)=(x\cdot yz)x^4$ (since cubes are in
the nucleus and $(xy)(zx)=(x\cdot yz)x$ is a Moufang identity).
We do not know whether an $M_4$ loop $L$ must be an A-loop.
By \cite{hala}, Theorem~1, $L$ is Moufang and $L/\Nuc(L)$ has
exponent three, but it is not clear whether $L/\Nuc(L)$ is
necessarily commutative.
We also do not know whether every loop isotope of a Moufang A-loop
is a (Moufang) A-loop. This would be true if
every $M_4$ loop is an A-loop, since 
the $M_k$ loops are isotopically invariant
(\cite{hala}, Theorem 2; \cite{pflugfelder}, IV.4.12)

Our investigations were aided by the automated deduction tool
OTTER developed by McCune \cite{mccune}; see Section \ref{sec:otter}
for further discussion.

\section{Preliminaries}

In preparation for the proofs of \lemref{lem:T-homom} and
\thmref{thm:main}, we now establish some notation and recall
some basic results from \cite{bp}. Let $L$ be a diassociative
A-loop. One can then derive many equations relating the
the $L(x,y)$, $R(x,y)$, and $T(x)$.

Define the permutation $J$ of $L$ by: $x J = x^{-1}$.
Conjugating by $J$, we have
$R(x)^{J}= J R(x) J = L(x^{-1})$;
likewise, $L(x)^{J}=R(x^{-1})$ and $L(x,y)^J = R(x^{-1}, y^{-1})$.
Note that $\varphi^{J}=\varphi$ for all 
automorphisms $\varphi$ of $L$; in particular for all
$\varphi \in \Mlt_1(L)$. Taking $\varphi = L(x,y)$, we have:
\begin{equation}\label{eq:L(x,y)=R(x',y')}
L(x,y)=R(x^{-1},y^{-1})
\end{equation}
for $x,y\in L$. Furthermore, from \cite{bp}, ((3.31) and (3.32)) we have
the following formulas for the inverses of the right and left
inner mappings:
\begin{align}
R(x,y)^{-1} &= R(y^{-1},x^{-1}) \label{eq:R(x,y)'=R(y',x')}\\
L(x,y)^{-1} &= L(y^{-1},x^{-1}) \label{eq:L(x,y)'=L(y',x')}
\end{align}

The fact that each $T(x)$ is an automorphism implies immediately:
\begin{align}
R(y) T(x) &= T(x) R(x^{-1} y x) \label{eq:RTTR} \\
L(y) T(x) &= T(x) L(x^{-1} y x) \label{eq:LTTL}
\end{align}

Another useful inner mapping is defined by
\begin{equation}\label{eq:C(x,y)}
C(x,y)=R(x)L(y)R(x^{-1})L(y^{-1}). 
\end{equation}
Since $C(x,y)^J = C(x,y)$, we also have:
\begin{equation}\label{eq:C-fixed}
C(x,y)=L(x^{-1})R(y^{-1})L(x)R(y).
\end{equation}
Also, by \cite{bp} (3.41):
\begin{equation}\label{eq:C-R}
C(x,y)=R(x,y)R(y,x)^{-1}.
\end{equation}

Further equations relating the $C(x,y),R(x,y),L(x,y)$ will 
be proved later (see Corollaries \ref{cor-clr} and \ref{cor-more-clr}).
As pointed out in \cite{bp}, in any loop, if $\varphi$ is an
automorphism which fixes an element $p$, then $\varphi$ commutes
with $L(p)$ and $R(p)$. In particular
(\cite{bp}, Lemma 3.3(i,ii,iii)), if $p,q,r$ are contained in
any subgroup of $L$, then:
\begin{align}
R(p)R(q,r) &= R(q,r)R(p)\text{; \ \ \ \ }L(p)R(q,r) = R(q,r)L(p)
\label{eq:com1}\\
R(p)L(q,r) &= L(q,r)R(p)\text{; \ \ \ \ }L(p)L(q,r) = L(q,r)L(p)
\label{eq:com2}\\
R(p)C(q,r) &= C(q,r)R(p)\text{; \ \ \ \ }L(p)C(q,r) = C(q,r)L(p).
\label{eq:com3}
\end{align}
One consequence is that the factors in the right and left inner
mappings can by cyclically permuted:
\begin{align}
R(x,y) &= R(y)R(y^{-1}x^{-1})R(x) = R(y^{-1}x^{-1})R(x)R(y)
\label{eq:permute1}\\
L(x,y) &= L(y)L(x^{-1}y^{-1})L(x) = L(x^{-1}y^{-1})L(x)L(y).
\label{eq:permute2}
\end{align}

\section{Proofs}

\begin{proof}[Proof of \lemref{lem:T-homom}]
For $x,y,z \in L$, we compute
\begin{align*}
z L(xy) T(x) &= y L(x) R(z) T(x) \\
&= y C(x^{-1},z^{-1}) R(z) L(x) T(x) && \by{eq:C-fixed} \\
&= y C(x^{-1},z^{-1}) R(z) R(x) \\
&= y R(z) R(x) C(x^{-1},z^{-1}) && \by{eq:com3} \\
&= y R(z) L(z^{-1}) R(x) L(z) && \by{eq:C(x,y)} \\
&= y T(z) R(x) L(z).
\end{align*}
By the mirror of this calculation and switching $x$ and $y$, we obtain:
\[
z R(xy) T(y^{-1}) = x T(z^{-1}) L(y) R(z) .
\]
But by (\ref{eq:RTTR}), we have
\[
y T(z) R(x) L(z) = y R(zxz^{-1}) R(z) = x T(z^{-1}) L(y) R(z).
\]
Hence, $L(xy) T(x) = R(xy) T(y^{-1})$, so that
$ T(x) T(y) = L(xy)^{-1} R(xy) = T(xy) $.
\end{proof}

\begin{corollary}
\label{cor-clr}
\begin{align}
R(x,y)^{-1} &= R(y,x) \label{eq:R(x,y)'=R(y,x)}\\
R(x,y) &= R(x^{-1},y^{-1}) = L(x,y) = L(x^{-1},y^{-1}) \label{eq:R=L}\\
C(x,y) &= C(x^{-1},y^{-1}) = R(x,y)^2. \label{eq:C(x,y)=C(x',y')}
\end{align}
\end{corollary}
\begin{proof}
(\ref{eq:R(x,y)'=R(y,x)}) follows from
\lemref{lem:T-homom} and (\ref{eq:bp(3.42)}).
To get (\ref{eq:R=L}), apply
(\ref{eq:R(x,y)'=R(y',x')}) and (\ref{eq:L(x,y)=R(x',y')}).
Then, (\ref{eq:C(x,y)=C(x',y')}) follows by using (\ref{eq:C-R}).
\end{proof}

\begin{lemma}\label{lem:key}
For all $x,y,z$ in a diassociative A-loop,
\begin{equation}\label{eq:key}
(yx)C(z,y)=(yx)C(z^{-1},x).
\end{equation}
\end{lemma}

\begin{proof}
Let $a = (yx)z^{-1}$. Then
\begin{align*}
(yx)C(z,y) &= (yx)C(z^{-1},y^{-1}) && \by{eq:C(x,y)=C(x',y')} \\
&= (yx)R(z^{-1}) L(y^{-1}) R(z) L(y) \\
&= a L(y^{-1}) R(z) L(y) \\
&= (y^{-1}a) R(a^{-1}(yx)) L(y) \\
&= (yx)L(a^{-1}) L(y^{-1}a) L(y) \\
&= (yx) L(y,a^{-1}) && \by{eq:permute2} \\
&= (yx) L(y^{-1},a) && \by{eq:R=L} \\
&= (ya^{-1})(ax) \\
&= (yx) R(x^{-1},a^{-1}) \\
&= (yx) L(x^{-1},a^{-1}) && \by{eq:R=L} \\
&= (yx) L(a^{-1}) L(xa) L(x^{-1}) && \by{eq:permute2} \\
&= x^{-1}(xa \cdot z) \\
&= (yx) R(z^{-1}) L(x) R(z) L(x^{-1}) \\
&= (yx)C(z^{-1},x).
\end{align*}
\end{proof}

\begin{proof}[Proof of \thmref{thm:main}]
For $x,y,z\in L$, we compute
\begin{align*}
x(y(xz)) &= xR(z)L(y)L(x)\\
&= xC(z,y)L(y)R(z)L(x)\\
&= (yx)C(z,y)R(z)L(x) && \by{eq:com3} \\
&= (yx)C(z^{-1},x)R(z)L(x) && \by{eq:key} \\
&= (yx)R(z)C(z^{-1},x)L(x) && \by{eq:com3} \\
&= (yx)L(x)R(z) \\
&= (xyx)z.
\end{align*}
\end{proof}

\begin{corollary}
\label{cor-more-clr}
$C(x,z) = L(z,x) = R(z,x)$, and $C(x,z)^3 = I$.
\end{corollary}
\begin{proof}
By the Moufang equation, $R(xz)L(x) = R(x)L(x)R(z)$.
Hence, $R(x^{-1}) R(xz) R(z^{-1}) = L(x)R(z) L(x^{-1}) R(z^{-1})$, so
that (by (\ref{eq:permute1}) and (\ref{eq:C-fixed}))
$R(z^{-1},x^{-1}) = C(x^{-1},z^{-1})$.
Now use Corollary \ref{cor-clr}.
\end{proof}

\section{Computer-aided Proofs}
\label{sec:otter}
We comment further on our use of McCunes program OTTER \cite{mccune}.
This is a general-purpose automated reasoning program which will
prove theorems from axioms in first-order logic.
In comparison with human reasoning, it is strongest in equational
reasoning, and weakest in domains such as set theory, where
there are many propositional connectives and alternations of quantifiers.
Thus, most of the \textit{new} mathematics
to come out of automated reasoning has been in fields
close to algebra. The book by Wos and Pieper \cite{wosp}
describes general methods for applying
automated reasoning to problems in mathematics and other areas.
Many new theorems proved by OTTER
occur in the book by McCune and Padmanabhan \cite{mcpad}.

Many authors (as in \cite{mcpad}) simply use the OTTER output
as the proof of a theorem. This is mathematically sound,
since although OTTER's search procedure is rather complex, 
the program can be made to output a simple \textit{proof object},
which can be independently verified by a short \verb+lisp+ program.
However, OTTER's proofs are often long sequences of
complicated equations which carry little intuitive content,
and it is useful to re-express them in a form which
a human reader can easily understand and verify.

Some discussion of the procedure for ``humanizing'' proofs
occurs in \cite{hk}. This was applied in the case of loop theory
in \cite{kuna, kunb, kunc, kund, kune}, and in the present paper,
where much of the argument is
cast in the spirit of Bruck and Paige \cite{bp}, emphasizing
group-theoretic properties of the $R(x)$ and $L(x)$,
rather than equations in the loop product and inverse.
For example, in Corollary \ref{cor-more-clr}, the statement
$C(x,z)^3 = I$ conveys more information to most human readers
than does the equivalent equation,
\[
z^{-1} (z ( ( z^{-1} (z ( ( z^{-1} (z ( y
x) x^{-1})) x) x^{-1})) x) x^{-1}) = y \ \ ,
\]
which might (in its \verb+ascii+ form) be a typical line of OTTER output.
However,
some proofs seem to require direct computations in the loop itself.
These proofs, although easy enough to verify by hand,
may lack some motivation.
The need for such computations probably explains why the results
of this paper have not been found before.

\begin{acknowledgements}
We wish to thank Toma\v{s} Kepka for suggesting this problem to us.
\end{acknowledgements}

\end{document}